\documentclass[11pt]{article}
\usepackage{epic,latexsym,amssymb}
\usepackage{color}
\usepackage{tikz}
\usepackage{amsfonts,epsf,amsmath}

\newtheorem{thm}{Theorem}
\newtheorem{lem}{Lemma}

\newcommand{\qed}{$\Box$}

\newcommand{\smallqed}{{\tiny ($\Box$)}}

\let\oldenumerate\enumerate
\renewcommand{\enumerate}{
  \oldenumerate
  \setlength{\itemsep}{0pt}
  \setlength{\parskip}{0pt}
  \setlength{\parsep}{0pt}
}

\begin{document}

\title{Niche Number of Linear Hypertrees}

\author{Thummarat Paklao, Nattakan Yahatta, Chutima Chaichana, \\Thiradet Jiarasuksakun and Pawaton Kaemawichanurat\thanks{This work was funded by Development and Promotion of Science Technology Talents (DPST) Research Grant 031/2559.}
\\ \\
Theoretical and Computational Science Center \\
Science Laboratory Building and Department of Mathematics\\
King Mongkut's University of Technology Thonburi \\
Bangkok, Thailand \\
\small \tt Email: pawaton.kae@kmutt.ac.th}

\date{}
\maketitle

\begin{abstract}
For a digraph $D$, the niche hypergraph $NH(D)$ of $D$ is the hypergraph having the same set of vertices as $D$ and the set of hyperedges is
\begin{align}
E(NH(D)) &= \{e \subseteq V(D) : |e| \geq 2~and~there~exists~v \in V(D)~such~that~e = N_{D}^{-}(v)\notag\\
         &~~~~~~~or~e = N_{D}^{+}(v)\}.\notag
\end{align}
A digraph is said to be acyclic if it has no directed cycle as a subdigraph. For a given hypergraph $H$, the niche number $\hat{n}(H)$ is the smallest integer such that $H$ together with $\hat{n}(H)$ isolated vertices is the niche hypergraph of an acyclic digraph. In this paper, we study the niche number of linear hypertrees with maximum degree two. By our result, we can conclude for a special case that if $H$ is a linear hypertree with $\Delta(H) = 2$ and anti-rank three, then $\hat{n}(H) = 0$. We also prove that the maximum degree condition is best possible. Moreover, it was proved that if $H$ is a hypergraph of rank $r$ whose niche number is not infinity, then $\Delta(H) \leq 2r$. In this paper, we give a construction of hypertrees whose niche number is $0$ of prescribed maximum degree from $3$ to $2r$.
\end{abstract}


\section{\bf Introduction}
A \emph{hypergraph} $H$ is an ordered pair $(V(H), E(H))$ where $V(H)$ is the set of elements and $E(H)$ is a family of subsets of $V(H)$. The elements of $V(H)$ are called the \emph{vertices} and the elements of $E(H)$ are called the \emph{hyperedges}. let $m = |E(H)|$ and $n = |V(H)|$. A \emph{loop} is a hyperedge containing exactly one vertex. A \emph{simple hypergraph} is a hypergraph such that if there exist two hyperedges $e, e'$ and $e \subseteq e'$, then $e = e'$. In this paper, all hypergraphs are simple and have no loop. A vertex which is not contained in any hyperedge is called an \emph{isolated vertex}. It is possible that our hypergraphs contain isolated vertices. The degree $deg_{H}(v)$ of a vertex $v$ in $H$ is the number of hyperedges that contain $v$. The maximum degree of $H$ is denoted by $\Delta(H)$. For a hyperedge $e$, the number of hyperedges of $H$ that intersect $e$ is called the \emph{edge degree} of $e$.
\vskip 5 pt

\indent The \emph{rank (anti-rank)} of a hypergraph is the maximum (minimum) number of vertices that each hyperedge can have. A hypergraph is said to be $r$-\emph{uniform} if both of rank and anti-rank are equal to $r$. The hypergraph is \emph{linear} if $|e \cap e'| \leqslant 1$ for all $e, e' \in E(H)$. For a subset $F$ of hyperedges of $H$, the hypergraph $H - F$ is obtained by removing all hyperedges of $F$ from $E(H)$ and removing all isolated vertices that might be occured. A hypergraph $H$ is \emph{connected} if, for any two vertices $u, v \in V(H)$, there exists a sequence of vertices $x_{1}, x_{2}, ..., x_{t}$ of $V(H)$ such that $u = x_{1}, v = x_{t}$ and, for $1 \leq i \leq t - 1$, we have $x_{i}, x_{i + 1} \in e$ for some $e \in E(H)$.
\vskip 5 pt

\indent A hypergraph is said to be a \emph{graph} if it is $2$-uniform. Hence, a graph $G$ is an ordered pair $(V(G), E(G))$ where $V(G)$ is the set of vertices and $E(G)$ is a set of subsets that has two vertices. A hyperedge of a graph is called an \emph{edge}. For a vertex subset $X$ of $V(G)$, the \emph{induced subgraph} $G[X]$ of $G$ on $X$ is the graph having the set of vertices $X$ and $xy \in E(G[X])$ if and only if $xy \in E(G)$. A \emph{tree} is a connected graph having no cycle as an induced subgraph.
\vskip 5 pt

\indent We now define a hypertree by using a tree. A hypergraph $H$ is a \emph{hypertree} if there is a tree $T$ such that $V(T) = V(H)$ and for each $e \in E(H)$ the induced subgraph $T[e]$ is connected. For a hypertree, a hyperedge having edge degree one is called a \emph{twig} and a hyperedge which is not a twig is called a \emph{trunk}. A vertex $v$ of $H$ is called a \emph{bud} if $deg_{H}(v) = 1$. For a linear hypertree $H$, if $e$ is a trunk of $H$ intersecting twigs $e_{1}, ..., e_{m}$, then the subhypergraph $\{e, e_{1}, ..., e_{m}\}$ of $H$ is called a \emph{branch}. It is possible that some branch has no twig. Remark that a linear hypertree is constructed from the intersection between branches.
\vskip 5 pt

\indent A graph which all the edges are oriented is called \emph{digraph}. An oriented edge of a digraph is called an \emph{arc}. All the arcs of digraphs in this work are oriented to be $(u, v)$ or $(v, u)$ but not both. If $(u, v)$ is an arc of a digraph, then there is an arrow from $u$ to $v$ where $u$ is the \emph{tail} and $v$ is the \emph{head} of the arrow. We let $A(D)$ be the set of all arcs of $D$. For a digraph $D$, we let $\overline{D}$ be the digraph which is obtained from $D$ by reversing all arcs of $D$. For a vertex $v$, the \emph{in-neighbour set} $N_{D}^{-}(v)$ of $v$ in $D$ is the set $\{ u : u \in V(D)$ and $(u, v) \in A(D)\}$ and the \emph{out-neighbour set}  $N_{D}^{+}(v)$ of $v$ in $D$ is the set $\{ u : u \in V(D)$ and $(v, u) \in A(D)\}$. The \emph{in-degree} of $v$ in $D$ is $|N_{D}^{-}(v)|$, similarly, the \emph{out-degree} of $v$ in $D$ is $|N_{D}^{+}(v)|$. For any two vertices $x_1$ and $x_p$, $p \geq 2$, of a digraph $D$, a \emph{directed path} $x_1,x_2,...,x_p$ is a sequence of different vertices such that $(x_1,x_2),(x_2,x_3),...,(x_{p-1},x_p) \in A(D)$. Moreover, a \emph{directed cycle} of $D$ is a sequence $x_1,x_2,...,x_p,x_1$ where $p \geq 3$ such that $x_1,x_2,...,x_p$ is a directed path of $D$ and $(x_p,x_1) \in A(D)$. An \emph{acyclic digraph} is a digraph having no directed cycle.
\vskip 5 pt

For an acyclic digraph $D$, the \emph{niche hypergraph} $NH(D)$ is the hypergraph with the vertex set $V(D)$ and the edge set
\begin{align}
E(NH(D)) &= \{e \subseteq V(D) : |e| \geq 2~and~there~exists~v \in V(D)~such~that~e = N_{D}^{-}(v)\notag\\
         &~~~~~~~or~e = N_{D}^{+}(v)\}.\notag
\end{align}

\noindent Let $I_{k}$ denote a set of $k$ isolated vertices. For a hypergraph $H$, the \emph{niche number} $\hat n(H)$ is the minimum number such that $H \cup I_{\hat n(H)}$ is the niche hypergraph of an acyclic digraph. If no such $k \in \mathbf{N}$ exists, we define $\hat n(H)= \infty$. Similarly, for an acyclic digraph $D$, the \emph{niche graph} $NG(D)$ is the graph with the vertex set $V(D)$ and the edge set
\begin{align}
E(NG(D)) &= \{xy : ~there~exists~v \in V(D)~such~that~(v, x), (v, y) \in A(D)\notag\\
         &~~~~~~~or~(x, v), (y, v) \in A(D)\}.\notag
\end{align}

\noindent For a graph $G$, the \emph{niche number} $\hat n_{g}(G)$ is the minimum number such that $G \cup I_{\hat n_{g}(G)}$ is the niche graph of an acyclic digraph.
\vskip 5 pt

\indent The niche graphs have been studied since 1989 by Cable et al. \cite{CJLS}. They gave some class of graphs whose niche number is infinity. They further established niche number of some special graphs such as a complete graph, a cycle and a path. They also proved that if $n_{g}(G) < \infty$, then $n_{g}(G) \leq |V(G)|$ for any graph G. In 1991, Bowser and Cable \cite{BC} further decreased the upper bound that if  $n_{g}(G) < \infty$, then $n_{g}(G) \leq \frac{2}{3}|V(G)|$. These works motivated Garske et al. \cite{CMH} to study niche hypergraphs. They constructed some hypergraphs whose niche number is infinity as given in the following. A \emph{hypernova} is a hypergraph $N(m)$ which $|E(N(m))| = m \geq 3$ and there exists exactly one vertex $x$ of $N(m)$ such that for $e, e' \in E(N(m))$, $e \cap e' = \{x\}$. They proved that :

\begin{thm}\cite{CMH}\label{thm 4}
If $N(m)$ is a hypernova, then $\hat n(N(m))=\infty$.
\end{thm}

\noindent Garske et al. \cite{CMH} further conjectured that the niche number of a hypercycle is $0$. Very recently, Kaemawichanurat and Jiarasuksakun\cite{KJ} proved that this conjecture is true.
\vskip 12 pt

\indent In this paper, we let
\vskip 5 pt

\indent $\mathcal{T} : $ the set of all linear hypertrees with maximum degree two and anti-rank three.
\vskip 5 pt

\noindent We prove that :

\begin{thm}\label{thm 1}
If $H \in \mathcal{T}$, then $\hat{n}(H) = 0$.
\end{thm}

The proof of Theorem \ref{thm 1} is given in Section \ref{proof}. In views of Theorem \ref{thm 4}, the hypernova $N(s)$ is a linear hypertree with $\Delta(N(s)) \geq 3$ and $\hat{n}(N(s)) = \infty$. Thus, the maximum degree condition in Theorem \ref{thm 1} is best possible. Although, there exists a linear hypertree with maximum degree at least three whose niche number is infinity, we can find a linear hypertree with maximum degree at least three whose niche number is $0$ too. Garske et al. \cite{CMH} proved that if $H$ is a hypergraph with $\hat{n}(H) < \infty$, then $\Delta(H) \leq 2r$ where $r$ is the rank of $H$. In Section \ref{proof2}, we prove the realizable problem that, for $r \geq 3$ and $3 \leq s  \leq 2r$, there exists a linear hypertree with maximum degree $s$ and rank $r$ whose niche number is $0$.
\vskip 5 pt

\begin{thm}\label{thm 2}
For $r \geq 3$ and $3 \leq s  \leq 2r$, there exists $r$-uniform linear hypertree $H$ with $\Delta(H) = s$ and $\hat{n}(H) = 0$
\end{thm}

\vskip 10 pt

\section{Proof of Theorem \ref{thm 1}}\label{proof}
In this section, we prove Theorem \ref{thm 1}. First of all, we may need the following lemma.
\vskip 5 pt

\begin{lem}\label{lem 0}
Let $H \in \mathcal{T}$. Then there exists an acyclic digraph $D$ which satisfies the following properties :
\begin{description}
  \item[(\emph{i})] $NH(D) = H$.
  \item[(\emph{ii})] For $e \in E(H)$, there is at most one bud $v \in e$ such that neither $N^{+}_{D}(v)$ nor $N^{-}_{D}(v)$ is an empty set.
\end{description}
\end{lem}

\noindent Remark that if $u \in e \setminus \{v\}$, then exactly one of $N^{+}_{D}(u)$ or $N^{-}_{D}(u)$ is an empty set. If $D$ is the acyclic digraph satisfying $(i)$ and $(ii)$ in Lemma \ref{lem 0}, then we say that $D$ is a \emph{good digraph of} $H$. It easy to see that Theorem \ref{thm 1} is a consequence of Lemma \ref{lem 0}. Thus, in this section, it suffices to prove Lemma \ref{lem 0}.
\vskip 12 pt

\noindent \emph{Proof of Lemma \ref{lem 0}.}\\
Because $\Delta(H) = 2$, it follows that $|E(H)| \geq 2$, in particular, $H$ is not $P_{1}$. We proceed the proof by induction on the number of branches of $H$. Hence, our basic step is when $H$ has exactly one branch. Therefore, $H$ has exactly one trunk $e$ that intersects some twig, $e_{1}, e_{2}, ..., e_{l}$ says. We may assume that $e = \{v_{1}, v_{2}, ..., v_{p}\}$ where $p \geq l$ and $e_{i} = \{v^{1}_{i}, ..., v^{r_{i}}_{i}\}$ such that $\{v^{1}_{i}\} = \{v_{i}\} = e \cap e_{i}$ for all $1 \leq i \leq l$.

\indent We now construct an acyclic digraph $D$ by setting $V(D) = V(H)$ and
\begin{center}
  $A(D) = \{(v^{2}_{1}, v) : v \in e\} \cup \{(v^{1}_{1}, v) : e_{l}\} \cup \{(v, v^{1}_{l}) : v \in e_{1}\} \cup (\cup^{l - 1}_{i = 2}\{(v, v^{r_{i - 1}}_{i - 1}) : v \in e_{i}\}).$
\end{center}

\noindent The digraph $D$ is illustrated in Figure 1 for example when $l = 4$.

\begin{figure}[h!]
\centering
\includegraphics[width=12cm]{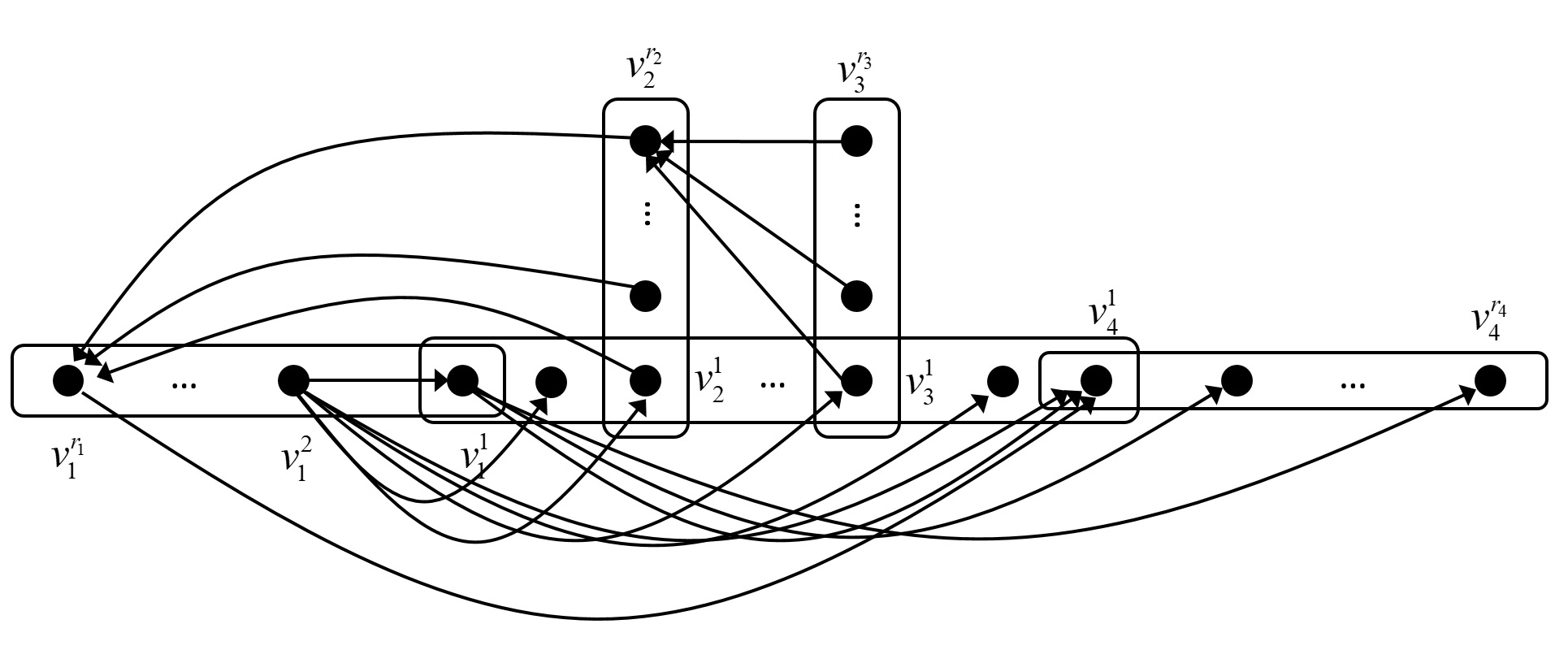}
\caption{\footnotesize{The digraph $D$ when $l = 4$.}}
\end{figure}

\noindent We see that $\{(v^{2}_{1}, v) : v \in e\}$ gives the hyperedge $e$, $\{(v^{1}_{1}, v) : e_{l}\}$ gives the hyperedge $e_{l}$, $\{(v, v^{1}_{l}) : v \in e_{1}\}$ gives the hyperedge $e_{1}$ and $\{(v, v^{r_{i - 1}}_{i - 1}) : v \in e_{i}\}$ gives the hyperedge $e_{i}$ for all $2 \leq i \leq l - 1$. It is not difficult to see that $D$ is an acyclic and giving no hyperedge apart form $e, e_{1}, ..., e_{l}$. Thus, $NH(D) = H$ and $D$ satisfies $(i)$.
\vskip 5 pt

\indent We will show that $D$ satisfies $(ii)$. If $v$ is a bud in $e$, then $N^{-}_{D}(v) = \{v^{2}_{1}\}$ and $N^{+}_{D}(v) = \emptyset$. We now consider $e_{i}$. When $i = 1$, we see that $N^{+}_{D}(v^{2}_{1}) = e$ and $N^{-}_{D}(v^{2}_{1}) = \emptyset$, $N^{+}_{D}(v^{j}_{1}) = \{v^{1}_{l}\}$ and $N^{-}_{D}(v^{j}_{1}) = \emptyset$ for all $3 \leq j \leq r_{1} - 1$, $N^{+}_{D}(v^{r_{1}}_{1}) = \{v^{1}_{l}\}$ and $N^{-}_{D}(v^{r_{1}}_{1}) = e_{2}$. Hence, the exactly one bud of $e_{1}$ whose in and out neighbors are not empty is $v^{r_{1}}_{1}$. When $i \in \{2, ..., l - 2\}$, we see that $N^{+}_{D}(v^{j}_{i}) = \{v^{r_{i - 1}}_{i - 1}\}$ and $N^{-}_{D}(v^{j}_{i}) = \emptyset$ for all $2 \leq j \leq r_{i} - 1$, and $N^{+}_{D}(v^{r_{i}}_{i}) = \{v^{r_{i - 1}}_{i - 1}\}$ and $N^{-}_{D}(v^{r_{i}}_{i}) = e_{i + 1}$. Hence, the exactly one bud of $e_{i}$ whose in and out neighbors are not empty is $v^{r_{i}}_{i}$. When $i \in \{l - 1, l\}$, we see that $N^{+}_{D}(v) = \{v^{r_{l - 2}}_{l - 2}\}$ and $N^{-}_{D}(v) = \emptyset$ for all buds $v \in e_{l - 1}$, and $N^{+}_{D}(v) = \emptyset$ and $N^{-}_{D}(v) = \{v^{1}_{1}\}$ for all buds $v \in e_{l}$. Thus, $D$ satisfies $(ii)$ and this establishes the basic step.
\vskip 5 pt

\indent Now, we may assume that every linear hypertree $H' \in \mathcal{T}$ with $1 \leq p' < p$ branches has a good digraph $D'$. Let $H \in \mathcal{T}$ having $p$ branches. We my need the following claim.

\noindent \textbf{Claim 1 :} There exists a branch $B$ such that $H - B$ is connected.\\
\indent Let $T$ be the set of all twigs of $H$ and consider $H - T$. Clearly, $H - T \in \mathcal{T}$ and all hyperedges of $H - T$ are the trunks of $H$. Let $e$ be a twig of $H - T$. Thus, $(H - T) - \{e\}$ is connected. Moreover, $e$ intersects some twigs of $H$ and $e$ together with these twigs form a branch in $H$, $B$ say. We see that $H - B$ is connected. This establishes Claim 1.\smallqed
\vskip 5 pt

\indent By Claim 1, we let $B$ be the branch that $H - B$ is connected. We further let $H' = H - B$. We may assume that the trunk of $B$ is $e$. Thus, $e$ intersects exactly one another trunk $f$ of $H'$. It is easy to see that $H'$ and $B$ are in $\mathcal{T}$, both of which have less than $p$ branches. By applying the inductive hypothesis to $H'$ and $B$, there exist good digraphs $D'$ and $D_{B}$ of $H'$ and $B$, respectively. Let $\{u\} = e \cap f$. Since $\Delta(H) = 2$, it follows that $deg_{B}(u) = deg_{H'}(u) = 1$, in particular, $u$ is a bud of $B$ and $H'$.
\vskip 5 pt

\indent In the following, we will show that there exist a good digraph $D''$ of $H'$ such that $N^{+}_{D''}(u)$ or $N^{-}_{D''}(u)$ is an empty set and a good digraph $D'_{B}$ of $B$ such that $N^{+}_{D'_{B}}(u)$ or $N^{-}_{D'_{B}}(u)$ is an empty set. After that, we construct an acyclic digraph $D$ of $H$ from the union of $D'_{B}$ and $D''$. We may establish the following claims.
\vskip 5 pt

\noindent \textbf{Claim 2 :} There exists a good digraph $D''$ of $H'$ such that $N^{+}_{D''}(u)$ or $N^{-}_{D''}(u)$ is an empty set.\\
\indent If $N^{+}_{D'}(u) = \emptyset$ or $N^{-}_{D'}(u) = \emptyset$, then we let $D'' = D'$ and this establishes the claim. Hence, we may assume that neither $N^{+}_{D'}(u)$ nor $N^{-}_{D'}(u)$ is empty set. By the property $(ii)$ of $D'$, $f$ has a vertex $x$ such that $N^{+}_{D'}(x) = \emptyset$ or $N^{-}_{D'}(x) = \emptyset$. Reverse all arcs of $D'$ if necessary, we may let $N^{+}_{D'}(x) = \emptyset$. We now let $D''$ be the digraph obtain from $D'$ by switching in and out neighbors between $u$ and $x$. That is $V(D'') = V(D')$ and
\begin{align}
  A(D'') &= (A(D') \setminus (\{(u, w), (w, u) : w \in N^{+}_{D'}(u) \cup N^{-}_{D'}(u)\} \cup \{(w, x) : w \in N^{-}_{D'}(x)\}))\notag\\
         &  ~~~~~ \cup (\{(x, w), (w, x) : w \in N^{+}_{D'}(u) \cup N^{-}_{D'}(u)\} \cup \{(w, u) : w \in N^{-}_{D'}(x)\}).\notag
\end{align}
Since we switch both in and out neighbors between $u$ and $x$ and $D'$ is acyclic, it follows that $D''$ is acyclic. It is not difficult to see that $NH(D'') = H'$ and $x$ is the exactly one bud of $f$ such that neither $N^{+}_{D''}(x)$ nor $N^{-}_{D''}(x)$ is empty set. Thus, $D''$ satisfies $(ii)$. This proves Claim 2.\smallqed
\vskip 5 pt

\noindent \textbf{Claim 3 :} There exists a good digraph $D'_{B}$ of $B$ such that $N^{+}_{D'_{B}}(u)$ or $N^{-}_{D'_{B}}(u)$ is an empty set.\\
\indent If $N^{+}_{D_{B}}(u) = \emptyset$ or $N^{-}_{D_{B}}(u) = \emptyset$, then we let $D'_{B} = D_{B}$ and this establishes the claim. Hence, we may assume that neither $N^{+}_{D_{B}}(u)$ nor $N^{-}_{D_{B}}(u)$ is empty set. By the property $(ii)$ of $D_{B}$, $e$ has a vertex $y$ such that $N^{+}_{D_{B}}(y) = \emptyset$ or $N^{-}_{D_{B}}(y) = \emptyset$. Reverse all arcs of $D_{B}$ if necessary, we may let $N^{+}_{D_{B}}(y) = \emptyset$. We now let $D'_{B}$ be the digraph obtain from $D_{B}$ by switching in and out neighbors between $u$ and $y$. That is $V(D'_{B}) = V(D_{B})$ and
\begin{align}
  A(D'_{B}) &= (A(D_{B}) \setminus (\{(u, w), (w, u) : w \in N^{+}_{D_{B}}(u) \cup N^{-}_{D_{B}}(u)\} \cup \{(w, y) : w \in N^{-}_{D_{B}}(y)\}))\notag\\
         &  ~~~~~ \cup (\{(y, w), (w, y) : w \in N^{+}_{D_{B}}(u) \cup N^{-}_{D_{B}}(u)\} \cup \{(w, u) : w \in N^{-}_{D_{B}}(y)\}).\notag
\end{align}
Since we switch both in and out neighbors between $u$ and $y$ and $D_{B}$ is acyclic, it follows that $D'_{B}$ is acyclic. It is not difficult to see that $NH(D'_{B}) = B$ and $y$ is the exactly one bud of $e$ such that neither $N^{+}_{D_{B}'}(y)$ nor $N^{-}_{D_{B}'}(y)$ is empty set. Thus, $D'_{B}$ satisfies $(ii)$. This proves Claim 3.\smallqed
\vskip 5 pt

\indent By Claims 2 and 3, we can always find a good digraph $D''$ of $H'$ such that $N^{+}_{D''}(u)$ or $N^{-}_{D''}(u)$ is an empty set and a good digraph $D'_{B}$ of $B$ such that $N^{+}_{D'_{B}}(u)$ or $N^{-}_{D'_{B}}(u)$ is an empty set. Reverse all arcs of $D''$ and $D'_{B}$ if necessary, we may assume that $N^{-}_{D'_{B}}(u) = \emptyset$ and $N^{+}_{D''}(u) = \emptyset$.
\vskip 5 pt

\indent Finally, we let $D$ be a digraph such that $V(D) = V(D'') \cup V(D'_{B})$ and $A(D) = A(D'') \cup A(D'_{B})$. Clearly, $NH(D) = H$. Thus, $D$ satisfies $(i)$. Moreover, $x$ is the exactly one bud of $f$ whose in and out neighbors are not empty while $y$ is the exactly one bud of $e$ whose in and out neighbors are not empty. By the Properties $(ii)$ of $D''$ and $D'_{B}$, we have that $D$ satisfies $(ii)$. Thus, $D$ is a good digraph of $H$. This completes the proof.
\qed

\section{Proof of Theorem \ref{thm 2}}\label{proof2}
In this section, we prove Theorem \ref{thm 2}. Recall that $N(m)$ is a hypernova which $|E(N(m))| = m$ and all hyperedges contain the vertex $a$. In this section, our $N(m)$ is $r$-uniform. Let $P_{5} = e_{1}, ..., e_{5}$ be an $r$-uniform linear hyperpath of $5$ hyperedges such that $e_{i} = \{v^{1}_{i}, v^{2}_{i}, ..., v^{r}_{i}\}$ and $v^{r}_{i} = v^{1}_{i + 1}$ for $1 \leq i \leq 4$. We also let $e' = \{u_{1}, u_{2}, ..., u_{r}\}$ be a hyperedge that contains $r$ vertices. An $(s - 1)$-\emph{petal flower} $F(s)$ is a hypergraph obtain from $N(s - 1), P_{5}$ and $e'$ by identifying the vertex $a$ of $N(s - 1)$ with the vertex $v^{1}_{1}$ of $e_{1}$ in $P_{5}$ and by identifying the vertex $u_{1}$ of $e'$ with the vertex $v^{2}_{2}$ of $e_{2}$ in $P_{5}$. We may call the hypernova part of $F(s)$ as the \emph{petals}. We may illustrate the 5-petal flower $F(6)$ by Figure 2.
\vskip 5 pt

\begin{figure}[h!]
\centering
\includegraphics[width=11cm]{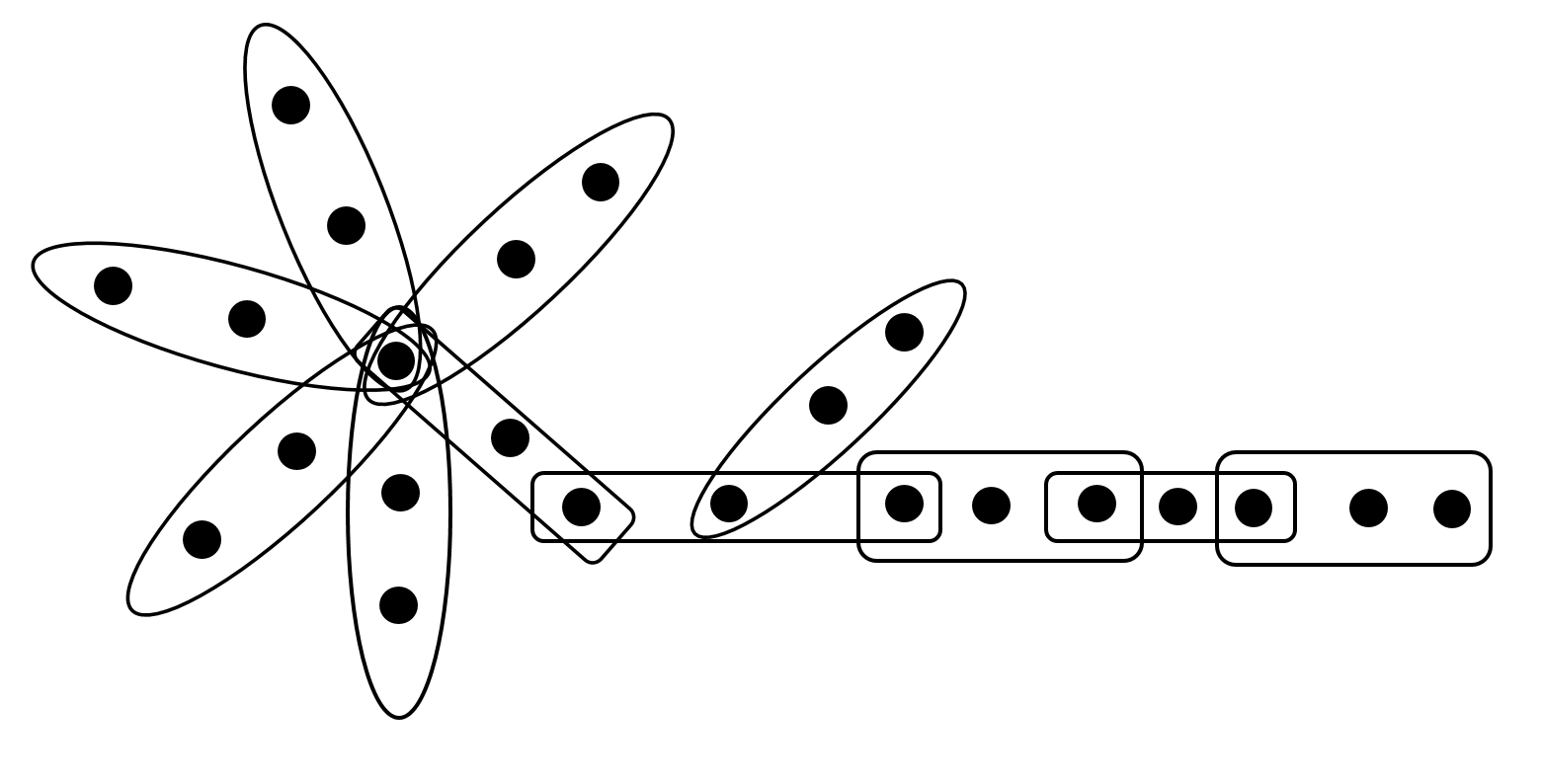}
\caption{\footnotesize{The 5-petal flower $F(6)$.}}
\end{figure}
\vskip 5 pt

\noindent Clearly, $deg_{F(s)}(a) = \Delta(F(s))= s$. We are ready to prove Theorem \ref{thm 2}.
\vskip 5 pt

\noindent \emph{Proof of Theorem \ref{thm 2}}.\\
We let $H = F(s)$. Clearly, $\Delta(H) = s$. We will show that $\hat{n}(F(s)) = 0$. We let $f_{1}, f_{2}, ..., f_{s - 1}$ be the petals of $F(s)$ such that $f_{i} = \{w^{1}_{i}, w^{2}_{i}, ..., w^{r}_{i}\}$ for $1 \leq i \leq s - 1$ and $w^{1}_{i} = a$. We may need the following claim.
\vskip 5 pt

\noindent \textbf{Claim 1 :} If $s \leq r$, then there exists acyclic digraph $D$ such that $NH(D) = F(s)$.\\
\indent Let $D$ be a digraph such that $V(D) = V(F(s))$ and
\begin{align}
  A(D) &= \{(v^{1}_{5}, v) : v \in e_{1}\} \cup \{(v, v^{1}_{4}) : v \in e_{2}\} \cup \{(v, v^{2}_{4}) : v \in e_{3}\} \cup \{(v^{1}_{3}, v) : v \in e_{4}\}\notag\\
         &  ~~~~~ \cup \{(v, v^{1}_{1}) : v \in e_{5}\} \cup \{(v^{1}_{1}, v) : v \in e'\} \cup (\cup^{s}_{i = 2}\{(v^{i}_{5}, v) : v \in f_{i - 1}\}).\notag
\end{align}
\noindent It is not difficult to check that $NH(D) = F(s)$. This establishes the claim.\smallqed
\vskip 5 pt

\indent By Claim $1$, we may assume that $r + 1 \leq s \leq 2r$, moreover, we let $D$ be an acyclic digraph such that $NH(D) = F(s) - f_{r} - f_{r + 1} \cdots - f_{s - 1}$. Thus, we let $D'$ be a digraph which $V(D') = V(D) \cup (\cup^{s - 1}_{i = r}f_{i})$ and
\begin{align}
  A(D') &= A(D) \cup (\cup^{s - r}_{i = 1}\{(v, u_{i}) : v \in f_{r - 1 + i}\}).\notag
\end{align}
\noindent It is not difficult to check that $NH(D') = F(s)$. This completes the proof.
\qed

\end{document}